\documentclass[number,citesort,dvips]{arxbj}
\usepackage{upgreek}

\aid{0}
\volume{17}
\issue{4}
\pubyear{2011}
\firstpage{1127}
\lastpage{1135}
\doi{10.3150/10-BEJ323}

\makeatletter
 \newtheorem{theoreme}{Theorem}[section]
\newtheorem{proposition}{Proposition}
\newtheorem{lemme}{Lemma}

\newremark{remarque}[theoreme]{Remark}

\newcommand{\R}{\mathbb{R}}
\newcommand{\1}{{\mathbf1}}
\newcommand{\pr}{\mathbb{P}}
\newcommand{\esp}{\mathbb{E}}
\newcommand{\N}{\mathbb{N}}
\newcommand{\cal}{\mathcal}
\newcommand{\F}{{\cal F}}

\makeatother

\begin{document}
\begin{frontmatter}

\title{First passage time law for some L\'{e}vy processes with compound
Poisson: Existence of a density}
\runtitle{First passage time law for some jump-diffusion processes}

\begin{aug}
\author[a]{\fnms{Laure} \snm{Coutin}\corref{}\thanksref{a}\ead[label=e1]{laure.coutin@math.univ-toulouse.fr}} \and
\author[b]{\fnms{Diana} \snm{Dorobantu}\thanksref{b}\ead[label=e2]{diana.dorobantu@univ-lyon1.fr}}
\runauthor{L. Coutin and D. Dorobantu}
\address[a]{IMT, University of Toulouse, Toulouse, France. \printead{e1}}
\address[b]{University of Lyon, University Lyon 1, ISFA, LSAF (EA
2429), Lyon, France.\\ \printead{e2}}
\end{aug}

\received{\smonth{11} \syear{2009}}
\revised{\smonth{7} \syear{2010}}

%
\begin{abstract}
Let $(X_t,  t\geq0)$ be a L\'{e}vy process with compound Poisson
process and $\tau_{x}$ be the first passage time
of a fixed level $x>0$ by $(X_t,  t\geq0)$. We prove that the law of
$\tau_x$ has a density
(defective when $\esp(X_1)<0)$ with respect to the Lebesgue measure.
\end{abstract}

%
\begin{keyword}
\kwd{first passage time law}
\kwd{jump process}
\kwd{L\'{e}vy process}
\end{keyword}

\end{frontmatter}

\section{Introduction}
\label{section1}
The main purpose of this paper is to show that the first passage time
distribution associated with a L\'{e}vy process with compound Poisson
process has a density with respect to the Lebesgue measure.

Let $X$ be a cadlag process started at $0$ and $\tau_x$ the first
passage time of level $x>0$ by~$X$.

L\'{e}vy, in \cite{Levy}, computed the law of $\tau_x$ when $X$ is a
Brownian motion with drift. This result is extended by Alili \textit{et
al.} \cite{alili} and Leblanc \cite{Leb} to the case where $X$ is an
Ornstein--Uhlenbeck process. The case where $X$ is a Bessel process was
studied by Borodin and Salminen in \cite{Bor}.

For the situation where the process $X$ has jumps,
the first results were obtained by Zolotarev \cite{zolotarev} and
Borokov \cite{borovkov} for $X$ a spectrally negative L\'{e}vy
process. Moreover, if $X_t$ has probability density $p(t,x)$ with
respect to the Lebesgue measure, then the law of $\tau_x$ has density
$f(t,x)$ with respect to the Lebesgue measure, where $xf(t,x)= tp(t,x)$
and $X_{\tau_x}=x$ almost surely.

If $X$ is a spectrally positive L\'{e}vy process, Doney \cite{doney}
gives an explicit formula for the joint Laplace transform of $\tau_x$
and the overshoot $X_{\tau_x} -x.$ When $X$ is a stable L\'{e}vy
process, Peskir \cite{Pes} and Bernyk \textit{et al.} \cite
{bernyk-dalang-peskir} obtain an explicit formula for the passage time density.

The case where $X$ has signed jumps has been studied more recently. In
\cite{doz}, the authors give the law of $\tau_x$ when $X$ is
the sum of a decreasing L\'{e}vy process and an independent compound
process with exponential jump sizes. This result is extended by Kou and
Wang in \cite{Kou} to the case of a diffusion process with jumps where
the jump sizes follow a double exponential law. They compute the
Laplace transform of $\tau_x$ and derive an expression for the density
of $\tau_x$. For a more general jump-diffusion process, Roynette \textit{et
al.} \cite{roy} show that the Laplace transform of $(\tau_x, x-X_{\tau
_{x_-}}, X_{\tau_x} -x)$ is the solution of some kind of random integral.

For a general L\'{e}vy processes, Doney and Kyprianou \cite
{doney-kyprianou} give the quintuple law of $(\bar{G}_{\tau_{x_-}},
\tau_x - \bar{G}_{\tau_{x}-}, X_{\tau_x} -x, x-X_{\tau_{x}-},
x-\bar{X}_{\tau_{x}-})$ where $\bar{X}_t = \sup_{s\leq t} X_s$ and
$\bar{G}_t =\sup\{s <t,  \bar{X}_s=X_s\}.$

Results are also available for some L\'{e}vy processes without Gaussian
component; see Lef\`{e}vre \textit{et al.} \cite{LL1,LL2,Pic1,Pic2}.
Blanchet \cite{Bl} considers a process satisfying the stochastic
equation $\mathrm{d}X_t=X_{t_-}(\mu \,\mathrm{d}t+\sigma\1_{\tilde{\phi}(t)=0}\,\mathrm{d}W_t+\phi
\1_{\tilde{\phi}(t)=\phi}\,\mathrm{d}\tilde N_t),  t\leq T,$ where $T$ is a
finite horizon, $\mu\in\R$, $\sigma>0$, $\tilde{\phi}(\cdot)$ is a
function taking two values, $0$ or $\phi$, $W$ is a Brownian motion,
$N$ is a Poisson process with intensity $\frac{1}{\phi^2}\1_{\tilde
{\phi}(t)=\phi}$ and $\tilde N$ is the compensated Poisson process.

The aim of this paper is to add to these results the law of a first
passage time by a L\'{e}vy process with compound Poisson process.

The paper is organized as follows: Section \ref{section2} contains the
main result (Theorem \ref{th1}) which gives the first passage time law
by a jump L\'{e}vy process. We compute the derivative of the
distribution function of $\tau_x$ at $t=0$ in Section \ref
{subsection1} and
at $t>0$ in Section \ref{subsection2}. Section \ref{appendice}
contains the proofs of some useful results.

\section{First passage time law}
\label{section2}

Let $m\in\R$ $(W_t,  t\geq0)$ be a standard Brownian motion $(N_t,
 t\geq0)$ be a Poisson process with constant positive intensity $a$
and $(Y_i,  i\in\N^*)$ be a sequence of independent identically
distributed random variables with distribution function $F_Y$ defined
on a probability space $(\Omega, {\mathcal F}, \pr).$ We suppose that
the $\sigma$-fields $\sigma(Y_i,  i\in\N^*)$, $\sigma(N_t,  t\geq
0)$ and $\sigma(W_t,   t\geq0)$ are independent. Let $(T_n,  n\in\N
^*)$ be the sequence of the jump times of the process $N$ and
$(S_i, i\in\N^*)$ be a sequence of independent identically
distributed random variables with exponential law of parameter $a$ such that
$T_n=\sum_{i=1}^n S_i$, $n \in\N^*.$

Let $\tilde{X}$ be the Brownian motion with drift $m\in\R$ and for
$z >0,$ $\tilde{\tau}_z=\inf\{t\geq0 \dvtx  mt+ W_t \geq z\}.$
By \cite{KS}, formula~(5.12), page 197, $\tilde\tau_z$ has the following
law on $\overline{\R}_+\dvtx
\tilde{f}(u,z)\,\mathrm{d}u+\pr(\tilde{\tau}_z=\infty)\delta_{\infty}(\mathrm{d}u),$
where\vspace*{-0.5pt}
%
\begin{eqnarray}
\label{derivbar}
\tilde{f}(u,z)&=&\frac{\mid z\mid}{\sqrt{2\uppi u^3}}\exp \biggl[-\frac
{(z-mu)^2}{2u} \biggr]\1_{]0, \infty[}(u), \qquad   u\in\R,  \quad
\mbox{and}\nonumber
\\[-8pt]
\\[-8pt]
\pr(\tilde{\tau}_z=\infty)&=&1-\mathrm{e}^{mz-|mz|}.
\nonumber
\end{eqnarray}
The function $\tilde{f}(\cdot,z)$ and all its derivatives admit $0$ as
right limit at $0$ and are $\cal{C}^{\infty}$ on $\R.$

Let $X$ be the process defined by
$ X_t=mt+ W_t+\sum_{i=1}^{N_t} Y_i,   t\geq0,$
and $\tau_{x}$ be the first passage time of level $x>0$ by $X\dvtx  \tau
_{x} =\inf\{u>0 \dvtx  X_u \geq x\}$.
The main result of this paper is the following theorem. 

\begin{theoreme}
\label{th1}
The distribution function of $\tau_{x}$ has a right derivative at $0$
and is differentiable at every point of $]0,  \infty[$. The
derivative, denoted $f(\cdot, x)$, is equal to
\[
f(0,x)=\frac{a}{2} \bigl(2 - F_Y(x)-F_Y(x_-) \bigr) +\frac
{a}{4}\bigl (F_Y(x)-F_Y(x_-) \bigr)
\]
and for every $t>0,$
\[
   f(t,x)=a\esp \bigl(\1_{\{\tau_{x}>t\}
}(1-F_Y) ( x-X_t ) \bigr) +\esp \bigl(\1_{\{\tau
_{x}>T_{N_t}\}}\tilde{f} (t-T_{N_t}, x-X_{T_{N_t}} ) \bigr).
\]
Furthermore, $\pr(\tau_{x}=\infty)=0$ if and only if $m+a\esp
(Y_1)\geq0.$
\end{theoreme}

The proof of Theorem \ref{th1} is given in Sections \ref{subsection1}
and \ref{subsection2}.


Let $(\F_t)_{t\geq0}$ be the completed natural filtration generated
by the processes $(W_t,  t\geq0)$, $(N_t,  t\geq0)$ and the random
variables $(Y_i,  i\in\N^*)\dvtx  \F_t=\sigma(W_s,  s\leq t)\vee
\sigma(N_s,  s\leq t,  Y_1, \ldots , Y_{N_t})\vee{\mathcal N}.$ Here,
${\mathcal N}$ is the set of negligible sets of $({\mathcal F},\pr).$

\begin{remarque}
This result is already known when $X$ has no positive jumps (see \cite{sato}, Theorem~46.4, page 348),
when $X$ is a stable L\'{e}vy process with no negative jumps (see \cite
{bernyk-dalang-peskir}) and
when $X$ is a jump diffusion where the jump sizes follow a double
exponential law (see \cite{Kou}).

According to \cite{LL2} and \cite{volpi}, for all $x>0,$
the passage time $\tau_x$ is finite almost surely if and only if
$m+ a \esp(Y_1) \geq 0.$
\end{remarque}

\subsection{Existence of the right derivative at $t=0$}
\label{subsection1}

In this section, we show that the distribution function of $\tau_{x}$
has a right derivative at $0$ and we compute this derivative.
For this purpose, we split the probability $\pr(\tau_{x}\leq h)$
according to the values of $N_h\dvtx
\pr(\tau_{x}\leq h)=\pr(\tau_{x}\leq h,  N_h=0)+\pr(\tau_{x}\leq
h,  N_h=1)+\pr(\tau_{x}\leq h,  N_h\geq2).$

Note that $\pr(\tau_{x}\leq h,  N_h\geq2)\leq1-\mathrm{e}^{-ah}-ah\mathrm{e}^{-ah}$
and thus $\lim_{h \rightarrow0}\frac{\pr(\tau_{x}\leq h,  N_h\geq2)}{h}=0.$

It suffices to prove the following two properties:
%
\begin{eqnarray}
\label{lem1}
 \frac{\pr(\tau_{x}\leq h,  N_h=0 )}{h}&\stackrel{h\rightarrow
0}\longrightarrow&0;\\
\label{lem2}
 \frac{\pr(\tau_{x}\leq h,  N_h=1 )}{h}&\stackrel{h\rightarrow
0}\longrightarrow&\frac{a}{2} \bigl(2-F_Y(x)-F_Y(x_-) \bigr) +\frac
{a}{4} \bigl(F_Y(x)-F_Y(x_-) \bigr).
\end{eqnarray}

On the set $\{\omega\dvtx  N_h(\omega)=0\}$, the processes $(X_t,   0 \leq
t \leq h)$ and $(\tilde{X}_t,   0 \leq t \leq h)$ are equal and $\pr
$-a.s. $\tau_{x}\wedge h =\tilde{\tau}_{x}\wedge h$. Since $\tilde
{\tau}_{x}$ is independent of $N$, we have
$\pr(\tau_{x}\leq h,  N_h=0 )=\mathrm{e}^{-ah}\pr(\tilde{\tau}_{x}\leq h).$
The law of $\tilde{\tau}_{x}$ has a $C^{\infty} $ density (possibly
defective) with respect to the Lebesgue measure, null on $]$--$\infty,
0]$, Thus, (\ref{lem1}) holds.

To prove (\ref{lem2}), we use the same type of arguments as in \cite
{roy} (for the proof of Theorem 2.4).
We split the probability $\pr(\tau_x \leq h, N_h=1)$ into three parts
according to the relative positions of $\tau_x$ and $T_1,$
the first jump time of the Poisson process $N$:
\begin{eqnarray*}
\pr( \tau_x \leq h,  N_h=1) &=& \pr( \tau_x <T_1,  N_h=1) + \pr
(\tau_x=T_1,  N_h=1)+ \pr(T_1<\tau_x \leq h,  N_h=1)\\
&=&A_1(h)+A_2(h)+A_3(h).
\end{eqnarray*}

\textit{Step} 1: As for (\ref{lem1}), we easily prove that $\frac
{A_1(h)}{h}\stackrel{h\rightarrow0}\longrightarrow0.$

\textit{Step} 2: We prove that $\frac{A_2(h)}{h}\stackrel{h\rightarrow
0}\longrightarrow\frac{a}{2}( 2 - F_Y(x) -F_Y(x_-)).$

Note that
$A_2(h)=\pr( \tilde{\tau}_x >T_1,  \tilde{X}_{T_1} + Y_1 \geq x,
 T_1 \leq h <T_2).$
Using the independence of $(S_i,  i\geq1)$ and $(Y_1$, $\tilde{X}$,
$\tilde{\tau}_{x}),$ we get
$\pr(\tau_{x}=T_1,  N_h=1)=a\mathrm{e}^{-ah}\int_0^h \esp (\1_{\{\tilde
{\tau}_{x}> s\}}\1_{\{Y_1\geq x-\tilde{X}_{s}\}} )\,\mathrm{d}s.$

Integrating with respect to $Y_1$, we obtain
\[
\frac{\pr(\tau_{x}=T_1,  N_h=1)}{a\mathrm{e}^{-ah}}=\int_0^h \esp\bigl (
(1-F_Y)\bigl((x-\tilde{X}_{s})_-\bigr) \bigr)\,\mathrm{d}s-\int_0^h \esp\bigl (\1_{\{
\tilde{\tau}_{x}\leq s\}}(1-F_Y)\bigl((x-\tilde{X}_{s})_-\bigr) \bigr)\,\mathrm{d}s.
\]

On the one hand, since $F_Y$ is a cadlag bounded function and $\tilde
{X}_s=ms+W_s$, where $W$ is continuous and symmetric, we get
$\lim_{s\rightarrow0}\esp (F_Y((x-\tilde{X}_{s})_-)
)=\frac{F_Y(x)+F_Y(x_-)}{2}.$
On the other hand, $\lim_{s\rightarrow0}\esp (\1_{\{\tilde{\tau
}_{x}\leq s\}}(1-F_Y)((x-\tilde{X}_{s})_-) )= 0.$

We deduce that
$\lim_{h \rightarrow0}\frac{A_2(h)}{h}=\frac{a}{2}
(2-F_Y(x)-F_Y(x_-) ).$

\textit{Step} 3: We prove that $\frac{A_3(h)}{h}\stackrel{h\rightarrow
0}\longrightarrow\frac{a}{4}( F_Y(x) -F_Y(x_-)).$

Note that
$\pr(T_1<\tau_{x}\leq h,  N_h=1)=\pr(T_1<\tau_{x}\leq h,  T_1\leq h<T_2)$
and $T_2=T_1+S_2 \circ\theta_{T_1}$, where $\theta$ is the
translation operator.

Moreover, on $\{T_1<\tau_{x}\leq h<T_2\}$, $X_s=X_{T_1}+\tilde
{X}_{s-T_1}\circ\theta_{T_1}$, where $T_1<s\leq h$ and $\tau
_{x}=T_1+\tilde{\tau}_{x-X_{T_1}} \circ\theta_{T_1}$. The
strong Markov property gives, with $\esp^{T_1}(\cdot)$ standing for $\esp
(\cdot\mid\F_{T_1})$,
\begin{eqnarray*}
A_3(h)&=&\esp \bigl(\1_{\{\tau_{x}> T_1\}}\1_{\{h\geq T_1\}}\esp
^{T_1} \bigl(\1_{\{\tilde{\tau}_{x -X_{T_1}}\leq h-T_1\}}\1_{\{h-T_1<
S_2\}} \bigr) \bigr)\\
&=& \esp \bigl(\1_{\{\tau_{x}> T_1\}}\1_{\{h\geq T_1\}
}\mathrm{e}^{-a(h-T_1)}\esp^{T_1}\bigl (\1_{\{\tilde{\tau}_{x-X_{T_1}}\leq
h-T_1\}} \bigr) \bigr)\\
&=&-\esp\bigl (\1_{\{\tilde{\tau}_{x}\leq T_1\leq h\}}\1_{\{
X_{T_1}<x\}}\mathrm{e}^{-a(h-T_1)}\esp^{T_1} \bigl(\1_{\{\tilde{\tau
}_{x-X_{T_1}}\leq h-T_1\}} \bigr) \bigr)\\
&&{}
+\esp \bigl(\1_{\{h\geq T_1\}}\1_{\{X_{T_1}<x\}} \mathrm{e}^{-a(h-T_1)}\esp
^{T_1}\bigl (\1_{\{\tilde{\tau}_{x-X_{T_1}}\leq h-T_1\}}
\bigr) \bigr).
\end{eqnarray*}
Since the distribution function of $\tilde{\tau}_{x}$ has a null
derivative at 0, we have
\[
\lim_{h\rightarrow0}\frac{1}{h}\esp \bigl(\1_{\{\tilde{\tau
}_{x}\leq T_1\leq h\}}\1_{\{X_{T_1}<x\}}\mathrm{e}^{-a(h-T_1)}\esp^{T_1}
\bigl(\1_{\{\tilde{\tau}_{x-X_{T_1}}\leq h-T_1\}} \bigr) \bigr)=0.
\]
It remains to show that
$\lim_{h \downarrow0} \frac{G(h)}{h}= \frac{a}{4} [ F(x)-F(x^-)],$
where
\[
G(h)=\esp \bigl(\1_{\{h\geq T_1\}}\1_{\{X_{T_1}<x\}}
\mathrm{e}^{-a(h-T_1)}\esp^{T_1} \bigl(\1_{\{\tilde{\tau}_{x-X_{T_1}}\leq
h-T_1\}} \bigr) \bigr).
\]

Integrating with respect to $T_1$ and then using the fact that $\tilde
{f}(\cdot, z)$ is the derivative of the distribution function of $\tilde
{\tau}_z$, we get
$G(h)=a\mathrm{e}^{-ah}\int_0^h\int_0^{h-s}\esp [\1_{\{\tilde
{X}_s+Y_1<x\}}\tilde{f}(u, x-\tilde{X}_s-Y_1) ]\,\mathrm{d}u\,\mathrm{d}s.$

We may apply Lemma \ref{lemme1_p2} to $p=1,$ $\mu=x-ms-Y_1$ and
$\sigma=\sqrt{s}.$ Then,
\[
\esp\bigl[\tilde{f}(u, \mu+\sigma G)\1_{\{\mu+\sigma G>0\}}\bigr]=\frac
{1}{\sqrt{2\uppi}}\esp \biggl[\mathrm{e}^{-{(\mu-mu)^2}/(2(\sigma
^2+u))} \biggl(\frac{\mu+\sigma^2m}{(\sigma^2+u)^{3/2}}+\frac
{\sigma G}{\sqrt{u}(\sigma^2+u)} \biggr)^+ \biggr]
\]
with $x^+=\max\{0, x\}$ and $G$ is a Gaussian $\cal{N}(0, 1)$ variable
and we have
\[
G(h)=\frac{a\mathrm{e}^{-ah}}{\sqrt{2\uppi}}\int_0^h\!\!\int_0^{h-s}\esp
\biggl[\mathrm{e}^{- {(x-m(u+s)-Y_1)^2}/(2(u+s))} \biggl(\frac{x-Y_1}{
(u+s)^{3/2}}+\frac{G\sqrt{s}}{\sqrt{u}(u+s)} \biggr)^+ \biggr]\,\mathrm{d}u\,\mathrm{d}s.
\]

We make the changes of variables $s=th,$ $u=hv$. Then,
\[
\frac{G(h)}{h}= \frac{a\mathrm{e}^{-ah}}{\sqrt{2\uppi}}\int_0^1\!\! \int_0^{1-t}
\esp \biggl[\mathrm{e}^{- {(x-mh(v+t)-Y_1)^2}/(2h(v+t))}
\biggl (\frac{x-Y_1}{ \sqrt{h}(v+t)^{3/2}}+\frac{G\sqrt{t}}{\sqrt
{v}(v+t)} \biggr)^+ \biggr]\,\mathrm{d}t\,\mathrm{d}v.
\]

However,
\[
\lim_{h \rightarrow0^+} \mathrm{e}^{-{(x-mh(T=v)-Y_1)^2}/(2h(t+v))}
\biggl(\frac{x-Y_1}{\sqrt{h}(t+v)^{3/2} }+\frac{G\sqrt{t}}{\sqrt
{v}(t+v)} \biggr)^+= \frac{\sqrt{t}}{\sqrt{v}(t+v)} G^+ {\mathbf
1}_{\{x=Y_1\}}
\]
and
\begin{eqnarray*}
&&\sup_{ 0 \leq h \leq1}\mathrm{e}^{-{(x-mh(t+v)-Y_1)^2}/(2h(t+v))}
\biggl(\frac{x-Y_1}{\sqrt{h}(t+v)^{3/2} }+\frac{G\sqrt{v}}{\sqrt
{1-v}} \biggr)^+\\
&& \quad  \leq\frac{\sup_{ z \geq0} z\mathrm{e}^{- {z^2}/{2}}
+|m|}{\sqrt{t+v}}+\frac{\sqrt{t}}{\sqrt{v}(t+v)}|G|.
\end{eqnarray*}

From Lebesgue's dominated convergence theorem,
we then obtain
\[
\lim_{h\rightarrow0}\frac{G(h)}{h}
=\Delta F_Y(x)\frac{\esp(G_+)}{\sqrt{2\uppi}}\int_0^1 \!\!\int
_0^{1-t}\frac{\sqrt{t}}{\sqrt{v}(t+v)}\,\mathrm{d}v\,\mathrm{d}t= \frac{1}{4}\Delta F_Y(x),
\]
where $\Delta F_Y(z)=F_Y(z)-F_Y(z_-).$ This identity achieves the proof
of step 3.

\subsection{Existence of the derivative at $t>0$}
\label{subsection2}

Our task now is to show that the distribution function of $\tau_{x}$
is differentiable on $\R_+^*$ and to compute its derivative.
For this purpose we split the probability $\pr(t<\tau_{x}\leq t+h),$
according to the values of $N_{t+h}-N_t$, into three parts:
\begin{eqnarray*}
&&\pr(t<\tau_{x}\leq t+h,  N_{t+h}-N_t=0)+\pr(t<\tau_{x}\leq t+h,
 N_{t+h}-N_t=1)\\
&& \qquad {}+\pr(t<\tau_{x}\leq t+h,  N_{t+h}-N_t\geq2)\\
&& \quad =B_1(h)+B_2(h)+B_3(h).
\end{eqnarray*}
Since $B_3(h)\leq\pr(N_{t+h}-N_t\geq2),$ we have $\lim_{h \rightarrow
0}\frac{B_3(h)}{h}=0$.

By the Markov property at $t$,
$B_2(h)=\esp(\1_{\{\tau_{x}> t\}}\pr^t(\tau_{x-X_t}\leq h,  N_h=1))$,
where $\pr^t(\cdot)=\pr(\cdot|\F_t).$

By (\ref{lem2}), $\frac{B_2(h)}{h}$ converges to
$\frac{a}{2}[2 -F_Y (x-X_t )-F_Y ( (x-X_t
)_- )]+\frac{a}{4}[F_Y (x-X_t )-F_Y (
(x-X_t )_- )]$ and is upper bounded by $\frac{\pr
(N_h=1)}{h}=a\mathrm{e}^{-ah}\leq a$. The dominated convergence theorem gives
\[
\lim_{h \rightarrow0}\frac{B_2(h)}{h}= a\esp \bigl(\1_{\{\tau_{x}>t\}
}(1-F_Y)(x-X_t) \bigr)+\frac{3a}{4}\esp \bigl(\1_{\{\tau_{x}>t\}
}\Delta F_Y(x-X_t) \bigr).
\]

However, the jumps set of $F_Y$ is countable and $X$ has a density (see
\cite{cont}, Proposition 3.12, page~90). Thus, $\esp (\1_{\{\tau
_{x}>t\}}\Delta F_Y(x-X_t) )=0$
and
$\lim_{h \rightarrow0}\frac{B_2(h)}{h}=a\esp (\1_{\{\tau_{x}>t\}
}(1-F_Y) (x-X_t ) ).$

It thus remain to prove that
%
\begin{equation}
\label{p1}
\frac{B_1(h)}{h}\stackrel{h\rightarrow0}{\longrightarrow}\esp
\bigl(\1_{\{\tau_{x}>T_{N_t}\}}\tilde{f} (t-T_{N_t},
x-X_{T_{N_t}} ) \bigr).
\end{equation}

Since $T_{N_t}$ is not a stopping time, we cannot apply the strong
Markov property. We split
\[
B_1(h)=
\pr(t<\tilde{\tau}_{x}\leq t+h<T_1)+\sum_{k=1}^{\infty}\pr
(t<\tau_{x}\leq t+h,  T_k< t<t+h<T_{k+1} ).
\]

On the set $\{T_k<t\}$, we have $X_t=X_{T_k}+X_{t-T_k}\circ\theta
_{T_k}$, hence on the set $\{\tau_x > T_k\},$ we have
$\tau_{x}=T_k+\tau_{x-X_{T_k}}\circ\theta_{T_k}.$ Moreover, on the
set $\{T_k< \min(t, \tau_x)\},$
\[
{\mathbf1}_{\{t<\tau_{x}\leq t+h,  T_k< t<t+h<T_{k+1}\}}= {\mathbf
1}_{\{ T_k<t\}}{\mathbf1}_{\{t-T_k<\tilde{\tau}_x\leq t+h -T_k
<S_{k+1} \}} \circ\theta_{T_k}
\]
and the
strong Markov property at $T_k$ gives
\begin{eqnarray*}
B_1(h)&=&\mathrm{e}^{-a(t+h)}\pr(t<\tilde{\tau}_{x}\leq t+h)\\
&&{}+\sum
_{k=1}^{\infty}\esp \bigl(\1_{\{T_k< t\}}\1_{\{\tau_{x}>T_k\}
}\mathrm{e}^{-a(t+h-T_k)}\esp^{T_k} \bigl(\1_{\{t-T_k<\tilde{\tau
}_{x-X_{T_k}}\leq t+h-T_k\}} \bigr) \bigr).
\end{eqnarray*}

The $\F_{T_k}$-conditional law of $\tilde{\tau}_{x-X_{T_k}}$
has the density (possibly defective) $\tilde{f}(\cdot, x-X_{T_k})$, thus
since $\mathrm{e}^{-a(t-T_k)}=\esp^{T_k} (\1_{\{T_{k+1}> t\}} )$, we have
\begin{eqnarray}
\label{eq-4}
B_1(h)&=&\mathrm{e}^{-ah}\int_t^{t+h}\esp \bigl( {\mathbf1}_{\{ 0 \leq t <T_1\}
} \bigr) \tilde{f}(u, x)\,\mathrm{d}u\nonumber\\
&&{}+\mathrm{e}^{-ah}\sum_{k=1}^{\infty}\int_t^{t+h} \esp \bigl( {\mathbf1}_{ \{
T_k \leq t <T_{k+1}\}}\1_{\{\tau_{x}>T_k\}}\tilde
{f}(u-T_k,x-X_{T_k}) \bigr)\,\mathrm{d}u
\\
&=& \mathrm{e}^{-ah}\int_t^{t+h}
\esp\bigl ( {\mathbf1}_{ \{T_{N_t} < \tau_x\}}\tilde{ f}(u-
T_{N_t}, x- X_{T_{N_t}})  \bigr)\,\mathrm{d}u .\nonumber
\nonumber
\end{eqnarray}
Since $\tilde{f}$ is continuous with respect to $u,$ for all $t>0,$
almost surely,
\[
\lim_{h \downarrow0} \frac{1}{h} \int_t^{t+h}{\mathbf1}_{\{
T_{N_t} < \tau_x\}}\tilde{ f}(u- T_{N_t}, x- X_{T_{N_t}})\,\mathrm{d}u= {\mathbf
1}_{\{ T_{N_t} < \tau_x\}}\tilde{ f}(t- T_{N_t}, x- X_{T_{N_t}}).
\]
According to Proposition \ref{cor8_p2} in the \hyperref[appm]{Appendix}, the family of
random variables $(\frac{1}{h} \int_t^{t+h}\tilde{ f}(u- T_{N-t}, x-
X_{T_{N_t}})\,\mathrm{d}u)_{0<h\leq1}$ is uniformly integrable.
We then obtain
\[
\lim_{h \rightarrow0}\frac{B_1(h)}{h}
=\esp\bigl (\1_{\{\tau_{x}>T_{N_t}\}}\tilde{f}(t-T_{N_t},
x-X_{T_{N_t}}) \bigr).
\]
Using (\ref{p1}), we deduce that
\[
\frac{\pr(t<\tau_{x}\leq t+h)}{h} \stackrel{h\rightarrow
0}\longrightarrow a\esp \bigl(\1_{\{\tau_{x}>t\}}(1-F_Y)(
x-X_t) \bigr)+\esp \bigl(\1_{\{\tau_{x}>T_{N_t}\}}\tilde
{f}(t-T_{N_t}, x-X_{T_{N_t}}) \bigr).
\]
The proof of Theorem \ref{th1} is thus complete.

\begin{appendix}\label{appm}
\section*{Appendix}\label{appendice}

We prove the following on $\tilde{f}$ given in (\ref{derivbar}).
\begin{lemme}
\label{lemme1_p2}
Let $G$ be a Gaussian random variable $\cal{N}(0,1)$ and let $\mu\in
\R$, $\sigma\in\R_+,$ $p \geq1$ and
$x^+=\max\{x,0\}.$ Then, for every $u\in\R,$
\begin{eqnarray*}
&&\esp\bigl[\tilde{f}(u, \mu+\sigma G)^p\1_{\{\mu+\sigma G>0\}}\bigr]\\
&& \quad =\frac
{1}{\sqrt{2^p\uppi^p}} \frac{u^{({1-2p})/{2}}\mathrm{e}^{-{p(\mu
-mu)^2}/(2(p\sigma^2+u))}}{ (p\sigma^2+u)^{ ({p+1})/{2}}} \\
&& \qquad {}\times\esp
\Biggl[ \Biggl(\sigma G + \sqrt{\frac{u}{p\sigma^2 +u}} (\mu- mu) + m\sqrt
{u(p\sigma^2 +u)} \Biggr)_+^p \Biggr].
\end{eqnarray*}
\end{lemme}


\begin{proposition}
\label{cor8_p2}
For every $t>0$ and $1 \leq p < 3/2,$
%
\[
\sup_{0 <h \leq1} \esp \biggl[  \biggl(\frac{1}{h}\int_t^{t+h}
{\mathbf1}_{ \{T_{N_t} < \tau_x\}}\tilde{ f}(u- T_{N_t}, x-
X_{T_{N_t}})\,\mathrm{d}u  \biggr)^p \biggr]< + \infty.
\]
\end{proposition}
\begin{pf}
Let $I(h)$ be
\[
I(h)=\frac{1}{h}\int_t^{t+h}
{\mathbf1}_{ \{T_{N_t} < \tau_x\}}\tilde{ f}(u- T_{N_t}, x-
X_{T_{N_t}})\,\mathrm{d}u.
\]
Using Jensen's inequality, the following estimate holds:
\[
\esp(I(h)^p ) \leq\frac{1}{h}\int_t^{t+h} \esp \bigl(
{\mathbf1}_{ \{x- X_{T_{N_t}}>0\}}\tilde{ f}(u- T_{N_t}, x-
X_{T_{N_t}})^p \bigr)\,\mathrm{d}u.
\]
Conditioning by the filtration generated by $N$ and $Y_i$, $i\in
{\mathbf N},$ it becomes, where $G$ is a standard Gaussian random
variable independent of $N$ and $Y_i$, $i\in{\mathbf N},$
\begin{eqnarray*}
&&\esp(I(h)^p ) \leq\frac{1}{h}\int_t^{t+h} \esp\Biggl (
{\mathbf1}_{ \{x- mT_{N_t} -\sum_{i=1}^{N_t} Y_i -\sqrt{T_{N_t}}
G>0\}}\\
&&\hphantom{\esp(I(h)^p ) \leq\frac{1}{h}\int_t^{t+h} \esp\Biggl (}{}\times\tilde{ f}\Biggl(u- T_{N_t}, x- mT_{N_t} -\sum_{i=1}^{N_t} Y_i -\sqrt
{T_{N_t}} G\Biggr)^p \Biggr)\,\mathrm{d}u.
\end{eqnarray*}

Note that for $u \in[t,t+h],$ $t -T_{N_t}\leq u-T_{N_t}\leq1+t
-T_{N_t},$ $pT_{N_t}+ t-T_{N_t} >t$ and if $C_{p}= \sup_{ x \in
{\mathbf R}^+} \sqrt{x^p} \mathrm{e}^{-{px}/{2}},$ then, from Lemma
\ref{lemme1_p2},
\begin{eqnarray*}
&&\esp(I(h)^p ) \leq\frac{3^{p-1}}{ \sqrt{2^p \uppi^p}}\esp \biggl(
\frac{ T_{N_t}^{ {p}/{2}}}{(t-T_{N_t})^{p- {1}/{2}}t^{
({p+1})/{2}}}\esp(|G|^p) + \frac{1}{(t-T_{N_t})^{ ({p-1})/{2}}t^{
{1}/{2}+p}}C_{p}\\
&&\hphantom{\esp(I(h)^p ) \leq\frac{3^{p-1}}{ \sqrt{2^p \uppi^p}}\esp \biggl(}{} + |m|^p \frac{1}{t^{ {1}/{2}}(t-T_{N_t})^{ ({p-1})/{2}}}
 \biggr).
\end{eqnarray*}
Observe that for every $t>0$ and $(\alpha,  \gamma)\in\,]$--$1, 0]\times
[0,+\infty[$, the random variables $(t-T_{N_t})^{\alpha}
T_{N_t}^{\gamma} $ are integrable (see the details below), which
completes the proof of Proposition \ref{cor8_p2}.

Note that
%
\begin{equation}\label{app-1}
\esp \bigl((t-T_{N_t})^{\alpha} T_{N_t}^{\gamma} \bigr)
\leq t^{\alpha} + \sum_{i=1}^{\infty}\esp \bigl(\1
_{t>T_i}(t-T_i)^{\alpha}T_i^{\gamma} \bigr)<+\infty.
\end{equation}
However, for $i \geq1,$ $T_i$ admits as density the function $u
\mapsto\frac{a^i}{(i-1)!}u^{i-1} \mathrm{e}^{-au},$ thus
\begin{eqnarray*}
\esp\bigl (\1_{\{t>T_i\}}(t-T_i)^{\alpha}T_i^{\gamma} \bigr)&=&\frac
{a^i}{(i-1)!}\int_0^t\mathrm{e}^{-au}(t-u)^{\alpha}u^{\gamma+i-1}\,\mathrm{d}u
\leq\frac{a^i}{(i-1)!}\int_0^t(t-u)^{\alpha}u^{\gamma+i-1}\,\mathrm{d}u\\
&=&\frac{a^i}{(i-1)!}t^{\gamma+i+\alpha}\frac{\Gamma(\gamma
+i)\Gamma(\alpha+1)}{\Gamma(\gamma+i+\alpha+1)} .
\end{eqnarray*}

Consequently, the sum in the right-hand term of inequality (\ref
{app-1}) is finite and the random variable $(t-T_{N_t})^{\alpha}
T_{N_t}^{\gamma} $ is integrable.
\end{pf}
\end{appendix}

\textbf{Acknowledgements}
The authors would like to thank M. Pontier
and P. Carmona for their careful reading, G. Letac for his helpful
comments and the referee for his helpful suggestions concerning the
presentation of this paper.

\printhistory

\end{document}